\newcommand{\R}{\mathbb{R}}
\newcommand\cB{\mathcal{B}}
\newcommand\cC{\mathcal{C}}
\newcommand\cD{\mathcal{D}}
\newcommand\cE{\mathcal{E}}

\newcommand\cL{\mathcal{L}}

\newcommand\cR{\mathcal{R}}

\newcommand\cT{\mathcal{T}}

\newcommand\ag{\alpha}

\newcommand\gd{\delta} 

\newcommand\go{\omega}
\newcommand\gO{\Omega}

\newcommand\gt{\theta}

\newcommand{\Hom}{\operatorname{Hom}}

\newcommand{\Id}{\mathrm{Id\,}}
\newcommand{\ve}{\varepsilon}
\newcommand{\ga}{\gamma}

\newcommand{\la}{\lambda}
\newcommand{\Ker}{\mathrm{Ker\,}}
\newcommand{\Coker}{\mathrm{Coker\,}}

\newcommand \inv{^{-1}}

\newcommand \xx{{x_*}}
\newcommand \yy{{y_*}}
\newcommand \ttt{{t_*}}

\newcommand\im{\operatorname{Im}}

\newcommand{\sgn}{\operatorname{sgn}}

\newcommand{\hand}{\hbox{ \,and\, }}
\newcommand{\hif}{\hbox{ \,if\,} }
\newcommand{\hof}{\hbox{\,of\,} }

\newcommand{\hforall}{\hbox{ \,for\  all\, }}
\newcommand{\LL}{\Lambda}

\newcommand{\sk}{\vskip10pt}

\newcommand{\bt}{\begin{theorem}}
\newcommand{\et}{\end{theorem}}
\newcommand{\bc}{\begin{corollary}}
\newcommand{\ec}{\end{corollary}}
\newcommand{\bp}{\begin{proposition}}
\newcommand{\ep}{\end{proposition}}
\newcommand{\bit}{\begin{itemize}}
\newcommand{\eit}{\end{itemize}}
\newcommand{\bl}{\begin{lemma}}
\newcommand{\el}{\end{lemma}}
\newcommand{\br}{\begin{remark}}
\newcommand{\er}{\end{remark}}
\newcommand{\bd}{\begin{definition}}
\newcommand{\ed}{\end{definition}}
\newcommand{\be}{\begin{equation}}
\newcommand{\bee}{\begin{equation*}}
\newcommand{\eee}{\end{equation*}}
\newcommand{\ee}{\end{equation}}
\newcommand{\barr}{\begin{array}}
\newcommand{\earr}{\end{array}}
\newcommand{\ra}{\rightarrow}
\newcommand{\bal}{\begin{align}}
\newcommand{\eal}{\end{align}}
\newcommand{\dd}{\colon}
\newcommand{\Rmp}{\R^{m+p}}

\newcommand{\Rm}{\R^{m}}
\newcommand{\Rp}{\R^{p}}

\newcommand{\vs}{\vskip5pt}

\newcommand{\xst}{x_*}
\newcommand{\mst}{m_*}
\newcommand{\pst}{p_*}

\documentclass[12]{amsart}
\usepackage[latin1]{inputenc}
\usepackage[T1]{fontenc}
\usepackage{amsmath}
\usepackage{amsfonts}
\usepackage{amsmath}
  \usepackage{paralist}
  \usepackage{graphics} %% add this and next lines if pictures should be in esp format
  \usepackage{epsfig} %For pictures: screened artwork should be set up with an 85 or 100 line screen
\usepackage{graphicx}  \usepackage{epstopdf}%This is to transfer .eps figure to .pdf figure; please compile your paper using PDFLeTex or PDFTeXify.
 \usepackage[colorlinks=true]{hyperref}
 \hypersetup{urlcolor=blue, citecolor=red}

\newtheorem{theorem}{Theorem}[section]
\newtheorem{corollary}{Corollary}

\newtheorem{lemma}[theorem]{Lemma}
\newtheorem{proposition}{Proposition}

\theoremstyle{definition}
\newtheorem{definition}[theorem]{Definition}
\newtheorem{remark}{Remark}

\title[Criticality and crisis] {REMARKS ON CRITICALITY AND CRISIS IN PURE EXCHANGE ECONOMIES}

\email{jacobo.pejsachowicz@polito.it}
\thanks{J. Pejsachowicz  is supported by GNAMPA-INDAM}

\author{Jacobo Pejsachowicz}
\address{Dipartimento di Scienze Matematiche\\
Politecnico di Torino\\
Corso Duca degli Abruzzi, 24\\
10129 Torino\\
Italy}
\email{jacobo.pejsachowicz@polito.it} 
\begin{document}
\dedicatory{Dedicated to the memory of A. Granas}
%\date{\today}
\begin{abstract} In the framework of Balasko's theory of  Equilibrium Manifold in a Pure Exchange Economy, we introduce a rigorous notion of crisis transition and establish some infinitesimal criteria which allow us to detect an unavoidable crisis among the critical equilibria. The proof builds on the relationship between the mathematical notions of branching, envelopes, and the intrinsic derivative. 
\end{abstract}

\maketitle

% DEFINITIONS, REMARKS NOT NUMBERED

%\newtheorem{definition}[theorem]{Definition}
%\newtheorem{example}[theorem]{Example}

\tableofcontents
%\begin{abstract}

\section*{INTRODUCTION}  

\begin{quote} \textit{ Y yo pienso que, para la Ciencia, hay que temer como una enfermedad el ir buscando una justificaci\'on fuera de s\'i misma, as\'i como hay que temer en el hombre la pregunta del fin, del  \textquestiondown para  qu\'e ? de la vida. Ustedes saben cuanta filosof\'ia  desesperada se halla en el fondo de esta pregunta; y sin embargo se vive por el amor a la vida, por el amor a los hijos, por el amor a la humanidad. As\'i es para la ciencia: las teor\'ias valen por la luz interior que han dado a quien las cre\'o, valen por la luz que dan todav\'ia   a quien  las estudia, no importa que esta luz pueda derivar de una  pregunta del entendimiento puro o de una de la ciencia aplicada.
 El fin de la vida es la vida digna y el fin de la ciencia es la ciencia digna; mas el juicio de la dignidad sale s\'olo de nuestra conciencia, por lo cual est\'an igualmente lejanas de la verdad ambas f\'ormulas: la de la ciencia para la pr\'actica y la de la ciencia para la ciencia.}
 \begin{flushright}\textrm{Beppo Levi} \end{flushright}
\end{quote}

This article ultimately deals with the relationship between three mathematical objects: branching, envelopes, and the intrinsic derivative. This relationship has been studied in Calculus of Variations, Differential Geometry and Bifurcation Theory at least since the introduction of conjugate points by  G.G.J.Jacobi. However, our aim here is to explain the possible relevance of this relationship to mathematical economics and, more precisely, to an approach to Walras equilibria initiated by Yves Balasko in the seventies.   Since then there has been considerable progress in topological bifurcation theory,  \cite{Iz}, and we believe that some of the results and methods of this theory can be helpful in understanding and completing Balasko's work. 
\vs

 In the Arrow-Debreu model of a Pure Exchange Economy, regular economies play a fundamental role. Identifying the economy with the vector whose components are the initial endowments of all participants to the market, these are the regular values of the projection to the commodity space restricted to the graph of the  Walras correspondence.  
  To consider regular values make sense since  Balasko \cite{Ba75} showed that the graph of the Walras correspondence is a submanifold $E$ of the Cartesian product of the price simplex $S$ with the commodity space $\Omega,$ called  \emph{The Equilibrium Manifold}.   

More precisely: according to Balasko, an \emph{ equilibrium} is a pair $e=(p,\go)\in S\times\gO$ where  $\go\in \gO$ is the vector of initial endowments of the market traders (in short an \emph {economy}),  and $p\in S$  is the associated  \emph{equilibrium price,} namely a price vector $p$ at which the aggregate excess demand $z(p,\go)$ of the market vanishes.  Denoting  with $W(\go)$ the set of \emph{equilibrium prices} associated to a given economy $\go,$ the \emph{equilibrium manifold} $E$ is the graph of the multivalued correspondence $W.$   In other words, the equilibrium manifold is the set $E =\{(p,\go) \in S\times \gO |\  z(p,\go)=0 \},$ endowed with the smooth structure of a submanifold of $S\times \gO.$ The \emph{natural projection,} or \emph {Debreu map,} $\pi\colon E\to \Omega$ is the restriction to $E$ of the projection of $S\times \gO$ onto $\gO.$ 
\vs

    A regular point of a map between manifolds of the same dimension is a  point where the differential of the map is an isomorphism between the corresponding tangent spaces. Regular points of $\pi$  are called \emph{regular equilibria } in mathematical economics.  The complement of the open set $\cR$ of regular equilibria is the set $\cC$ of critical equilibria.  A \emph{regular economy} is an endowment $\go\in \gO$ such that every point in the inverse image  $ \pi ^{-1} (\go) $ is a regular equilibrium. The main virtue of regular economies is due to the fact that for every economy $\go'$ close enough to $\omega$ one can select a unique, smoothly depending equilibrium price $p(\go')$  close to a given equilibrium price $p$ associated with $\go.$  Indeed, under appropriate  assumptions on the utility functions  the map $\pi\colon E\to \Omega$ is proper, and this together with the invertibility of the differential  implies that, for every regular economy $\go,$ the  fiber $E_{\go} =\pi ^{-1} (\go)  $ is a finite subset of $ E,$  and moreover every element $(p,\go) \in E_{\go} $ has an open  neighbourhood in $E$ diffeomorphic to a fixed neighbourhood  $N$ of $\go$ in $\gO.$  The inverse diffeomorphism composed with the projection  to  $S$  defines a smooth local selection of equilibrium prices on $N.$ 

However, not the smoothness but rather the continuity of price selections that are considered by Balasko and others as fundamental for a well-behaved trade market. 

 \emph{"The discontinuity property of equilibrium price selections contradicts an assumption that underlies, implicitly at least, many policy-oriented fields as, for example, international trade theory or public economics"}\cite{Ba09}.
 
 Critical (or singular) economies, i.e, those which contain at least one critical equilibrium on their fiber are frequently tied, at least heuristically, to economical crises understood from the mathematical viewpoint as a sudden jump of equilibrium prices of the market due to the discontinuity of any local choice of a unique equilibrium price in a neighbourhood of an equilibrium belonging to  $E_{\go}.$

 \emph{"This apparently inexplicable and unpredictable discontinuity leads to the serious, sometimes heated questioning of the market mechanism, and even to irrational behaviour that occasionally ends in widespread destruction of resources through futile attempts to get back to the former price levels"}\cite{AcPu}. 

  \vs

However, there is a long-overdue need to clarify the above interpretation of critical equilibria in terms of the economical crises. To some extent, see Appendix A, our interest in this issue arose from a pictorial description in the pages of \cite{Ba88} of the "futility of the attempts to restore prices", during an imaginary crisis of the market of artichokes in Bretagne.

\vs

 Let us notice that our discussion of regular equilibria shows that the noninvertibility of $D\pi(p,\go)$  is only a necessary but not sufficient condition for the appearance of jumps in equilibrium prices.  In some cases, even when $D\pi(p,\go)$ is noninvertible one can still associate continuously to every economy in a small enough neighbourhood of $\go$ a unique equilibrium price close to $p$ because a smooth map can be a local homeomorphism even at a critical point.  One can easily build examples of critical economies with locally defined continuous selections of equilibrium prices (Appendix B).  As a matter of fact, the latter assertion holds for any isolated critical equilibrium if  $\dim \gO$ is greater than three \cite{Pl}.

In this note, taking into serious consideration the above heuristics, we deeper refine the analysis of critical equilibria by introducing the notion of crisis transition together with the mathematical apparatus which permits to detect among the critical equilibria the ones that lead to a crisis. 
%Our approach is mainly topological and uses only some elementary facts from singularity theory \cite{Ar}.
\vs 
We define  \emph{crisis transition} (or simply  \emph{crisis})  as a point   $e_*=(p_*,\go_*)\in E $
such that no continuous selection $p(\go)$ of prices with $p(\go_*)= p_*$  exist on any neighbourhood of $\go_*.$ This is clearly equivalent to  say that $\pi$ fails to be a local homeomorphism at the point $e_*,$ and hence   
 the set of crisis transitions is nothing but the set $\cB$ of \emph{branch points} of the map $\pi.$
  A \emph{crisis economy} is any economy $\omega\in \gO$  containing among the points of the fiber $E_{\go}$ a crisis.  Since not every critical equilibrium is a crisis, the problem arises in finding criteria that distinguish crisis transitions from other critical equilibria. Here we will obtain some sufficient conditions for a given critical equilibrium to be a crisis that depend only on derivatives of $\pi$ at $e_*$ up to the second order.

While the precise formulation of these infinitesimal conditions will be stated in the next section, roughly speaking they can be described as follows:  we will show that an equilibrium $e_*\in E$  is a crisis whenever the differential  $D\pi(e_*)$ has an odd-dimensional kernel and the intrinsic derivative of the differential $D\pi$ at $e_*,$ in some direction tangent to $E,$ is an isomorphism.

 The above two conditions are far from being a characterization of crisis transitions. Indeed, using multi-parameter bifurcation theory one can find more general sufficient conditions for a crisis.  However, our criterion is interesting on its own, because, as we will show,  any equilibrium $e_* =(p_*,\go_*)$  verifying the above conditions  is not only a crisis but an \emph{unavoidable crisis} in the following sense: 
 in every neighbourhood of  $e_*$ there exists two regular equilibria $e_1=(p_1,\go_1)$ and $e_2=(p_2,\go_2)$ such that every continuous path  in $E$ between $e_1$ and $e_2$ must traverse a crisis.  Equivalently, every path in $\gO$ joining the two close economies $\go_1$ and $\go_2$ that can be lifted to a path in $E$ between $e_1$ and $e_2$ must contain some crisis economy. To say it in other terms, our infinitesimal condition at the point $e_*$ implies that the set of crisis transitions $\cB$ separates $E$ and $e_i, i=1,2$ belong to different paths components of its complement. 

Let  $\cE =\pi(\cB)$ be the set of crisis economies.  The natural projection restricts to a proper map $\pi \colon E \setminus \pi^{-1}(\cE) \ra \gO\setminus \cE$ which is a local homeomorphism. It is well-known \cite{Ma} that such a map is a covering projection and hence possesses the unique path lifting property.   Assuming that $\go_i; i=1,2$ are regular, as a consequence of our result, we have the following alternative: any continuous path of economies between  $\go_1$ and $ \go_2$  either intersects the set $\cE$ of crisis economies or the unique lifting of this path with initial point $e_1=(p_1,\go_1)$ must have an endpoint $e=(p,\go_2)$ with the price $p$ far from   $p_2.$   The above accounts for the futility in trying to bring back the equilibrium price $p_1$ to a  price close to $p_2$ by continuously modifying endowments without traversing a crisis economy. This underscores the importance of out of equilibrium dynamics in price adjustment problems. 

Our second result deals with the opposite case. We show that whenever no unavoidable crisis occurs on $E$ and the set of crisis transitions $\cB$ is "thin" in a very precise sense, then every economy either has a unique equilibrium price or it is indeterminate in terms of comparative statics. 
 
Unavoidable crises also bear relation to the appearance of the local Leontief's transfer paradox i.e., the possibility of improving the utility of the equilibrium allocation of a given trader by donating to other participants a small amount of his endowment.  While the local transfer paradox can be formulated only for regular economies, our results together with those of \cite{To} allow us to conclude that every equilibrium verifying the two conditions of an unavoidable crisis can be arbitrarily approximated by unstable regular economies admitting a local transfer paradox. We leave a detailed discussion of this relation to the interested reader.

Our final result is the comparison between the approach to crisis transitions via Balasko's natural projection and a more classical viewpoint in terms of the critical points of the excess demand map. A side outcome of this comparison is a simple criterion for the appearance of an unavoidable crisis in terms of the aggregate excess demand.
\vs
The above is a rough description of the economic content of the paper and the main results. However, it would be unfair to relegate the mathematics of the manuscript to a mathematical appendix,  as the economists often do, because it constitutes the main body and interest of the article.  But, on the other hand, it would be incorrect to reduce the remaining part to an application to mathematical economics. This simply wasn't how we thought about the problem.  
  
Accordingly, we decided to divide the article into two equally ranked parts. The first is devoted to the relevance of the above mathematical objects to economic theory,  discussing our main contribution on this matter along the lines of Balasko's approach to Walras equilibria.  In this part, we keep a full allegiance to the notation of his book \cite{Ba16}, from now on quoted as FTGE in the text. 

In the second part, we introduce the mathematical apparatus needed for the proof of the results stated in the first part and prove them in a more abstract mathematical framework, keeping in mind eventual applications to other fields. Here, we use the standard notation in nonlinear analysis, in particular reversing the order in which the variables and parameters appear. We also avoided the symbols used in the first part to spare the reader of 
any further confusion.

Intending to make the content more accessible to a mixed public, we reduced the mathematical background to a minimum by working on finite dimensions only. Even so, some knowledge of manifolds and vector bundles appears to be unavoidable. Using the degree theory of \cite{FPR} and the methods of \cite{AP, FP}  most of our results can be extended to Fredholm maps between Banach manifolds with potential applications to the infinite horizon and infinite-dimensional economies \cite{Co, Ba97, Ac}. 

\vs 

Here are the details of both parts of the manuscript.  In the first part, we shortly review Balasko's approach to Walras equilibria. Then we introduce the intrinsic derivative and state the main results of the paper: sufficient conditions for an unavoidable crisis, the analysis of markets with thin avoidable sets of crisis transitions, and the relation between the intrinsic derivative of the natural projection and the intrinsic derivative of the reduced aggregate excess demand.  

We begin the second part with an Envelope Theorem which provides a  sufficient condition for points of the discriminant of a family of plane curves to belong to the envelope and relate it with bifurcation and branching. Next, we introduce the intrinsic derivative and prove the first two theorems of the first section in a more abstract mathematical setting. Finally, we apply a reduction property of the intrinsic derivative to implicitly defined manifolds in Euclidean spaces and state our criteria for unavoidable branching in terms of the defining map.

There are three appendices. Appendix A is a description, taken from \cite{Ba88}, of the anomalous behavior of market prices during a crisis. Appendix B contains examples of critical economies which admit a continuous selection of equilibrium prices and crisis economies. Appendix C is a short review of Brouwer's degree and Krasnoselkij's Bifurcation Principle. 
\vs\vs

\section{CRITICAL EQUILIBRIA  AND UNAVOIDABLE CRISES}  
  
 \subsection{The equilibrium manifold and the natural projection}   Differential topology methods were introduced in economics by Debreu \cite{De} and Smale \cite{Sm71} in the seventies under the influence of the spectacular advances of differential topology and global analysis of that time. The extra assumption of smoothness on utility functions needed in their approach had as a positive counterpart the local uniqueness of the associated price in the case of regular economies. Balasko \cite{Ba75} developed a closely related approach based on analysis of the equilibrium manifold and the natural projection.  He also introduced singularity theory methods and the catastrophe theory of Thom in mathematical economics  \cite{Ba78}.  Since the definitions adopted by Balasko and his followers partially disagree with known standards, e.g., Balasko's equilibrium manifold is different from the one of Smale \cite{Sm74}, we shortly review some of his basic definitions and results, providing in this form a dictionary to his lingo.     
\vs
 The Arrow-Debreu model of pure exchange economy with $m$ consumers and $l$ goods  associates to each consumer an \emph{allocation, or commodity bundle}, belonging to the \emph{consumption space} $X = \R^l_{++},$  the strictly positive orthant of $\R^l,$ and a \emph{ price vector} $p=(p_1,\dots, p_l)$ belonging to the positive orthant $\R^l_{++}.$    We will consider the \emph{normalized prices}.  Choosing as numeraire the $l$-th commodity, these are  $p_i = p'_i/p'_l$  which make the normalized price of the numeraire equal to one.  Under  this choice,  the  price  "simplex" becomes  the hyperplane $S =\{ (p_1,\dots,p_l)\in \R^l_{++}| \ p_l=1\}.$  
 
 The $ i$-th  consumer participating to the market is represented by his \emph {utility function}  $u_i \colon X\ra \R$ which is smooth and satisfies  the following assumptions:
\begin{itemize}
\item[i)]Monotonicity: the gradient  $\nabla u_i (x)$ of $u_i$ at every point $x$ belongs to $\R^l_{++} .$ 
\item[ii)] Strict quasi-concavity:  The Hessian matrix $D^2 u_i(x)$ is negative definite on each tangent hyperplane  to the indifference surface $u_i(x) =r.$ 
\item[iii)] Necessity: The indifference surface $u_i(x)=r$ is closed in $\R^l$ for every $r\in \R$.
\end{itemize} 

 The  consumer  enters the market with an \emph{ initial  endowment} $\omega_i=(\omega_i^1,\dots,\omega_i^l),$  and  his goal  is to reach through trade
 an  \emph {allocation} $x_i=(x^1_i,\dots, x^l_i)\in X $  which maximizes his utility $u_i$ subjected to his {\bf budget constraint}  \[B_i=\{ x\in X : \ \sum p_k x_i^k\leq  wi, \   \text {where} \ w_i=\sum  p_k \go_i^k  \}. \] 

 Under the above assumptions  the  utility function $u_i$ assumes  its maximum on  the boundary $\partial B_i=\{ x\in X| \ p\cdot x=w_i\},$ of the budget set   at a unique point $ x=x_i(p,\go_i)$  which represents the \emph{ demand } of the $i$-th consumer under  the price system $p$ and initial endowment $\go_i.$ 
  The  function $x_i\colon S\times \R \ra X$   is smooth and satisfies the Walras law  $p\cdot x_i(p,\go_i) = w_i.$ 
  
 Define \emph{ economy } $\go =(\go_1,\dots\go_m) \in \gO=\R^{ml}_{++},$ as 
the vector of initial endowments of all participants to the market. 

In this article, the utility functions of the consumers will be assumed to be fixed, while the economy $\go$  will be considered as an exogenous parameter of the market.  

\bd

The \emph{ excess demand} of the market is the map  $z\colon S\times \gO\ra X$ defined by  \be \label{excess}  z(p,\go_1,\dots\go_m) = \sum_{I=1}^m (x_i(p,\go_i)-\go_i).\ee 
\ed

As a consequence of the Walras law for individual demand functions also the excess demand  $z$ verifies the same property. Namely,  for every $\go \in \gO,$
\be \label{walras} p\cdot z(p,\go)=0.\ee 

\bd An \emph{equilibrium} is a pair $e=(p,\go)\in S\times\gO$ such that $z(p,\go)=0.$\ed
According to Balasko \cite{Ba09} $\go$ is the economy and $p$ is the equilibrium price associated with the equilibrium $e.$ The monotonicity, strict quasi-convexity, and necessity assumptions imply that the set E of Walras equilibria is a submanifold of dimension $ml$ of $S\times \gO,$ and that the natural projection $\pi $ associating to every equilibrium its economy is a proper surjective map. 

 The proof of this goes, roughly speaking, as follows: as a consequence of Walras law, the last coordinate of the map $z$ is a linear combination of the first $l-1$ coordinates. Therefore E is also the set of zeros of the reduced aggregate excess demand $\bar z$  obtained from $z$ by eliminating the last coordinate. A direct calculation shows that  $0$ is a regular value of the map $\bar z$ and hence $E= \bar z^{-1}(0)$ is a submanifold by the Implicit Function Theorem. That $\pi$ is proper follows from the necessity axiom and properties of individual demand functions. Finally, the surjectivity and hence existence of an equilibrium for every initial endowment follows from homotopy invariance of Brouwer's degree (see the Appendix) and the uniqueness of the equilibrium price associated with a non-trade economy, i.e., the one whose initial endowment is already the preferred allocation.  It can be shown that $E$ is diffeomorphic to the Euclidean space and in particular, is connected.

We already discussed some of  Balasko's terminology in the introduction and here we will stick to that one. Regular equilibria can be equivalently defined either as regular points  $e=(p,\omega)$ of  $\pi$ or points  whose first coordinate $p$ is a regular point of the partial map $\bar z_\go\colon S\ra \R^{l-1}$ \cite[Proposition 4.5.4.]{Ba09}.  Regular economies are the regular values of the natural projection. Using properness as in \cite{Mi} one shows that every regular economy $\go$  possesses an evenly covered neighbourhood, i.e., an open neighbourhood $N$ such that $\pi^{-1} (N)$ is a  finite, disjoint union of open subsets of $E$ diffeomorphic by  $\pi$ with $N.$  Composing the inverse with the projection on $S$ one gets a smooth local selection of prices close to any equilibrium price associated to $\go.$
 Non-regular equilibria and economies are called \emph{critical}. They are seldom considered in economics, however, they have been studied in comparative statics using singularity theory, mainly in connection with the change of the number of equilibria of regular economies \cite{APP, Ba78}. 

For completeness, we will add another useful property of regular economies. Namely, that the set  $\cR$ of regular economies is open and that every economy can be arbitrarily approximated by a regular one.  The first assertion holds because proper maps send closed sets to closed sets and hence the set of critical economies $\cC=\gO\setminus \cR,$  being the image of the set of critical equilibria is closed. The second is a consequence of  Sard's theorem according to which  $\cC$  is of Lebesgue measure zero.  
\vs

\subsection{Existence of unavoidable crises} 
\bd. A \emph{crisis transition}, or simply  a \emph{crisis}  is a point  $e_*=(p_*,\go_*)\in E$ such that no continuous selection of prices $p(\go)$  with   $p(\go_*)=p_*$ exist on any  neighbourhood  of  $\go_*.$\ed 

As we observed in the introduction the above is equivalent to saying that the natural projection $\pi \colon E\ra  \gO$ fails to be a local homeomorphism at $e_*. $ 
The latter condition is precisely the definition of a \emph{ branch point } of a map.  
Thus, the set $\cB$ of crisis transitions is nothing but the set of branch points of $\pi.$
It follows from the Implicit Function Theorem that the set $\cB $ of crisis transitions is a closed subset of the set $\cC =E \setminus \cR$ of critical equilibria. However,  in general, $\cB$ is only a proper subset and we want to find sufficient conditions for a given critical equilibrium to be a crisis. 
 
Our main result is the following theorem which, from assumptions involving only the first and second-order derivatives of $\pi$ at a single point of  $E$, allows us to conclude that the set  $\cB$  of crisis transitions separates $E.$

 \bt \label{t1}
 If $e=(p,\go)$ is an equilibrium such that 
 \bit 
 \item [i)] $\dim \Ker D\pi (e)$  is odd 
  \item [ii)] for some $v\in T_e(E),$  the second intrinsic derivative  $\partial^2_v \pi (e)$ in the direction $v$  is an isomorphism. 
  \eit
  
  Then $e$ is an unavoidable  crisis. Namely,  it is a crisis and   in every neighbourhood of  $e$ in $E$ there exists two (regular) equilibria $e_\pm=(p_\pm,\go_\pm)$ belonging to different path components  of $E\setminus \cB.$ Equivalently,  every continuous path in $E$ between $e_+$ and $e_-$ must intersect the set of crisis transitions $\cB.$ 
  
  \et \begin{proof}  The equilibrium manifold  $E$ is connected and indeed diffeomorphic to a Euclidean space by Proposition 5.17 of FTGE.  In particular, it is orientable.  By Proposition 6.2  of FTGE,  the natural projection $\pi \dd E\ra\gO$ is proper. Thus the hypotheses of Theorem \ref{t4}  are verified, and Theorem \ref{t1} follows from this theorem.\end{proof}

  The set  $\cE =\pi(\cB)$ of all crisis economies is closed because proper maps send closed sets into closed ones.   By restriction $\pi$  induces a map $\pi\colon  E\setminus \pi^{-1} (\cE) \ra \gO\setminus \cE$ which is clearly proper and a local homeomorphism. Proper local homeomorphisms are covering maps and hence they possess the unique lifting property.  From the conclusions of the previous theorem, we get the following alternative which accounts for the impossibility to bring back the equilibrium price $p$ to a close price $p'$  by continuously modifying endowments without traversing a crisis. 

\bc [Impossibility to restore the prices] Let $e_\pm=(p_\pm,\go_\pm)$  be the regular equilibria near $e_*$ given by the previous theorem such that $\go_\pm$ are regular economies. Let  $d = dist( p_+, \pi\inv(\go_+)\setminus \{p_+\}).$  Then, any path $\ga$ from $\go_-$ to $\go_+$ either intersects the set $\pi({\cB})$ of crisis economies or the unique lifting of this path with initial point $e_-=(p_-,\go_-)$ must have an endpoint  $e=(p,\go_+)$ with  $\|p-p_+\|\geq d$ \ec 
\includegraphics[width=0.7\textwidth]{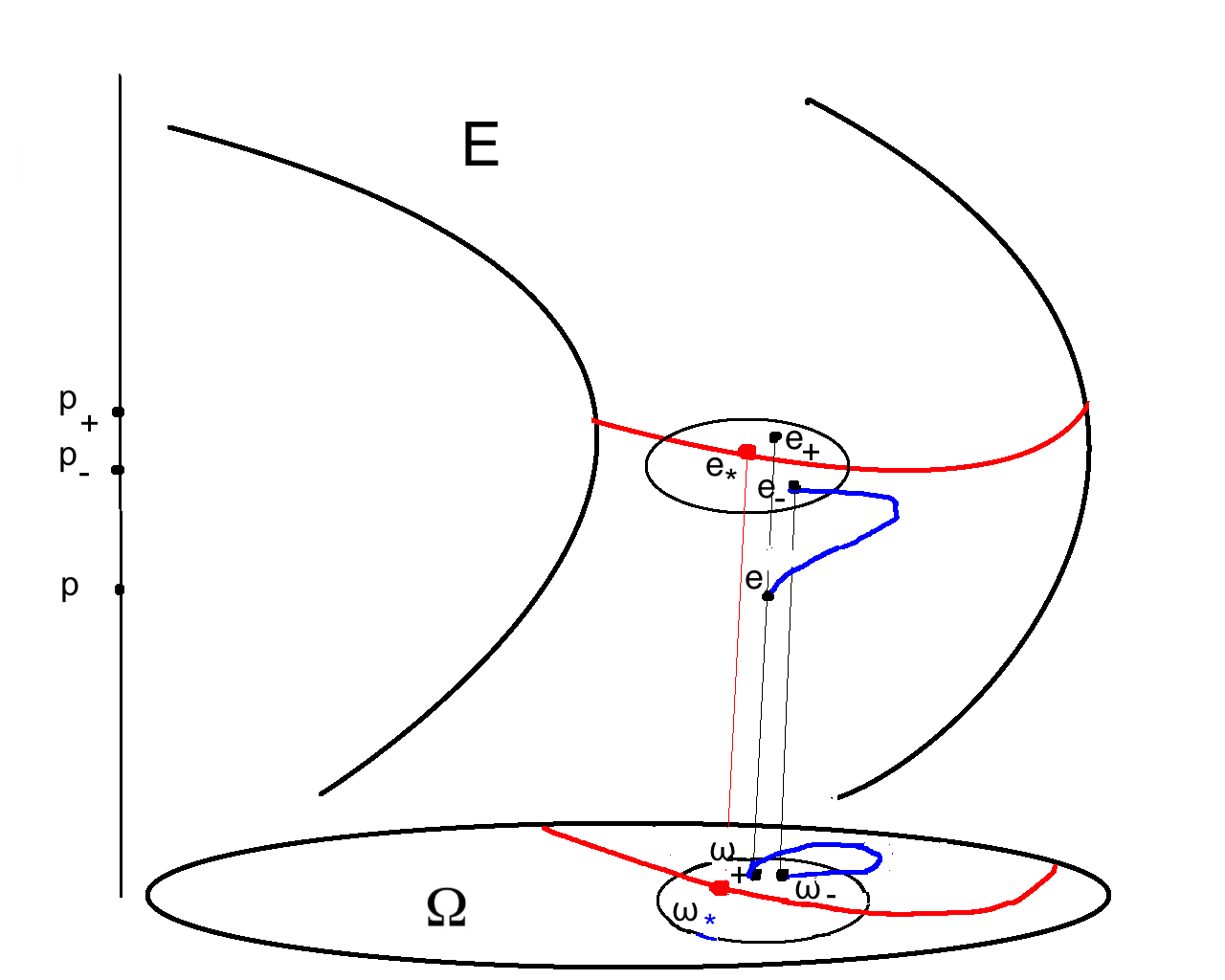}
\br Since the set of regular economies is dense in $\gO,$ any neighbourhood of $e_*$ in $E$ contains a pair  $e_\pm$ with the above property.\er
\proof.  Since $\pi \colon E\setminus \cB \ra \gO\setminus \pi(\cB)$ is a covering map, the path $\ga$ has a unique lifting  $\tilde \ga$ to $E$ with $\tilde \ga(0)= e_-.$ By Corollary \ref{c6} \, $\tilde \ga(1) \in \pi\inv(\go_+)\setminus \{p_+\}),$  which proves the assertion.
\vs

\subsection{Markets with thin sets of avoidable crisis transitions} 
Our next theorem explores what happens when the set of crisis transitions neither separates the manifold $E$ nor it contains interior points (closed sets with no interior points are said to be nowhere dense). It turns out that these two properties of $\cB$ lead to the uniqueness of prices of regular economies and the connectedness of the fibers over the crisis economies. 

\begin{theorem}\label{t18} Assume that the complement of the set of crisis transitions 
 $E \setminus \cB$ is both connected and dense in $E.$  Then for every $\go\in \gO,$ the fiber  $E_\go = \pi\inv (\go)$ is a continuum, i.e., a compact connected set.  In particular, if $E_\go$ contains some regular equilibrium $e$, then $e$ is the only element of this set. \end{theorem}  
 
 The mathematical property of local uniqueness of prices (or equivalently equilibria), typical of the regular economies, is termed \textit{ determinateness} in comparative statics.
Thus in terms of economic theory, the above theorem implies in particular that whenever the set of crisis transitions is nowhere dense and unavoidable branching does not occur,  then every economy either has a unique equilibrium price or is indeterminate.  
 
   \begin{proof} By  Proposition 6.34 of  FTGE, the Brouwer degree of the natural projection $\pi$ is $+1.$ Now the result follows from Theorem \ref{t8} of the second part. \end{proof}
 
 \vs 
\subsection{Excess demand and unavoidable crises}
 It is easy to see that a point $e= (p,\go)\in E $ is a critical point of the natural projection $\pi$ if and only if $p$ is a critical point of the reduced aggregate excess demand map 
$\bar z_\go\dd S\ra \R^{l-1},$ defined by $\bar z_\go(p):=\bar z(p,\go),$ and moreover the differentials at both points have canonically isomorphic kernels.  Our final result is the computation of the intrinsic second derivative $ \pi $ in terms of the second differential of $ \bar z_\go.$   This computation uses, together with a natural reduction property of the intrinsic derivative, some kind of duality which we tried to motivate in the second part, using envelope theory for plane curves, but without much success, probably because we are still far from fully grasping its real nature.  However, the result is of practical interest because it provides a simple criterion for the existence of an unavoidable crisis in terms of the first and second differential of $\bar z_\go$ at $p.$

\bt \label {t5} Given a  point $ e=(p,\go)\in E$  the inclusion $ i(v) := (v,0)$ sends isomorphically  $\Ker D \bar z_{\go}(p)$ into $\Ker D\pi(e).$ In particular, $e=(p,\go)$ is a critical point of $\pi$ if and only if $p$ is a critical point of $\bar z_{\go}.$

Moreover, there exists  a canonical isomorphism $j\dd \Coker D\bar z_{\go}(p) \ra\Coker D\pi(e)$ such that  for every $v\in \Ker D \bar z_{\go}(p),$\be\label{e0} \partial ^2_{i(v)} \pi(e ) =j \partial^2_v\bar z_{\go}(p). \ee \et

\begin{proof} The first part follows from: 
\be  \Ker D\pi(e) =\Ker \Pi \cap T_{e}(E) =\R^{l-1} \times \{0\}\cap \Ker Df(e)= i(\Ker D\bar z_{\go}(p)).\ee 
The second assertion is the thesis of the Theorem \ref{t20} in Part II.\end{proof} 

 \vs 

As a consequence, spelling out the definition of  $\partial^2_v\bar z_{\go}(p)= \partial_v D\bar z_{\go}(p),$ we obtain the following criterion:

\bt \label{t6}

  If  $e= (p, \go) \in E$ verifies
   \bit 
 \item [i)] $\Ker \bar z_{\go}(p)$  is odd-dimensional. 
  \item [ii)]  There exists a vector  $v\in \Ker \bar z_{\go}(p),$ such that 
 \[D^2\bar z_{\go}(p)[v,u] \notin \im  D\bar z_{\go}(p) \hforall u \in \Ker \bar z_{\go}(p), u\neq 0.\] 
 \eit
  Then $e$ is an unavoidable  crisis, i.e., the conclusions of Theorem {\rm \ref{t1}} hold.  
  \et
  
  \proof To obtain the result from  Theorem \ref{t21}  we have to verify the hypotheses of that theorem.  We have only to show that $E$ is locally bounded over $\gO.$ This property would an easy consequence of each $\bar z_\go$ being a proper map. Indeed, on finite dimensions,  properness of $\bar z_\go$ on $S$ is equivalent to  $\|\bar z_\go\|\to \infty,$ when $p$ goes either to infinite or to the boundary of $S.$ 
By continuity the above condition at $p$ implies the existence of a bound for $E$ on a small enough neighbourhood of $p.$  Unfortunately, as is shown in page 135 of FTGE, the map $\bar z_\go$ is not proper. However, using an idea of Balasko, we substitute $\bar z$ with the map $\tilde z$ defined by 
 \be \tilde z(p_1,\dots, p_{l-1},\go)=\left( (1+p_1) z^1(p_1,\dots, p_{l-1},\go),\dots,  (1+p_{l-1})z^{l-1}(p_1,\dots, p_{l-1},\go)\right)\ee     
 Then,  by Lemma $7.17$ of FTGE,  $\tilde z_\go $ is proper  for all $\go\in \gO$  and clearly $E=\tilde z\inv (0).$ Thus  $E$ is locally bounded over $\gO$ and the conclusion follows from Theorem \ref{t21}.

 \vs
\section{INTRINSIC DERIVATIVES AND BRANCHING}

\subsection{The envelope of a family of plane curves}

There are in the mathematical literature several notions of envelope of an implicitly  defined  one-parameter family  of  plane curves  $ C_z = \{ (x,y) | f(x,y,z) =0 \}.$ Assuming that  $f$ is  smooth 
the \emph{ envelope} of the family is  frequently defined in one of the following four forms:  
\bit 
\item The set $\cE$ of limit points of intersections of nearby curves in the family.\footnote{"intersection of two consecutive curves of the family " according to the picturesque language of Giuseppe Peano in \cite{Pe}.} 
\item The plane curve  $\cE'$ that is tangent to every curve of the family.
\item The boundary  $\cE''$ of the subset of the plane filled by all curves of the family. 
\item The set $\cD=\{(x,y) | f(x,y,z)=0, \frac{\partial f}{\partial z}(x,y,z)=0\}.$ \eit
In general   $\cE, \cE', \cE'', \cD$   differ one from the other \cite{BG}, but in some cases, they all coincide with $\cD.$  A well-known example in which this happens  is the envelope of ballistic trajectories:

 \includegraphics[width=0.8\textwidth]{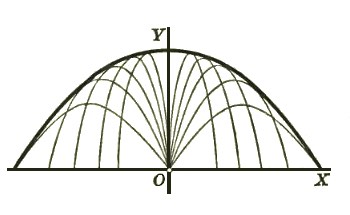}

\[f(x,y,z )=x\tan z -{\frac  {gx^{2}}{2v^{2}\cos ^{2}z }}-y=0.\]

 We will first shortly review the precise relationship between the first one $\cE,$ called here the \emph{envelope} of the family and the last one $\cD,$ which we will call the \emph{discriminant}.  As a partial motivation for the extension of this theory to higher dimensions, we will also consider the relation of envelopes with bifurcation and branching.
 
By definition a  point $q_*=(x_*,y_*)$ belongs to $\cE$  if there exists a  sequence of intersection points $q_n\in C_{z_n}\cap C_{z'_n}$ with $ z_n\neq z'_n $ but $ z_n, z'_n \ra z_* $ and $q_n \ra q_*$.  Notice that for every intersection point $q_n= (x_n,y_n)$ it holds that  
\be {\left(f(x_n,y_n, z_n) -f(x_n,y_n,z'_n)\right) /(z_n-z'_n) =0.}\ee If $z_* = \lim z_n,$ then  $f(x_*,y_*,z_*)=0$ by continuity. While the Lagrange intermediate value theorem allows us  to conclude that  $ \frac{\partial f}{\partial z}(x_*,y_*,z_*)=0.$ Thus $q_* \in \cD$  which proves that  $\cE\subset\cD.$

 A partial converse to the above is the following  theorem \cite{Pe} : 
\bt \label{t2}(The  envelope theorem)  Let $q_* \in \cD\cap C_{z_*}$ be such that, denoting the partial derivatives as subindices, the determinant    
\be \Delta(q_*) = [f_x f_{yz} - f_y f_{xz}](x_*,y_*,z_*) \neq 0, \ee
 then  $q_*\in \mathcal E.$  Moreover the curve $C_{z_*}$ is a regular curve in a neighbourhood of $q_*,$ and there exist a parametrization  $\ga\colon (z_*-\gd, z_*+\gd) \ra \R^2$ of  $C_{z_*}$  with $\ga(z_*) =q_*$  such that $\ga(z) $ is the unique  intersection point of $C_{z_*}$ with  $C_z.$
  \et  
  \proof    Define 
  \be g(x,y,z) =  \begin{cases}\frac{ \left(f(x,y, z) -f(x,y,z_*)\right)}{(z-z_*)} & \hif  z\neq z_* \\  \frac{\partial f}{\partial z}(x,y,z_*)  &\hif z=z_* .\end{cases}
  \ee
Consider the  map $F(x,y,z) =(f(x,y,z), g(x,y,z)).$  Clearly  $F(x_*,y_*,z_*)=(0,0)$ and  the Jacobian determinant of $F(-,-,z_*)$ at $q_*=(x_*,y_*)$ is $\Delta(x_*,y_*)\neq 0.$ By the Implicit Function Theorem there exists a parametrized curve $\ga(z)=(x(z),y(z))$ defined  in $(z_*-\gd, z_*+\gd)$ such that $\ga(z_*)=q_*$ and $F(\ga(z),z)=(0,0).$ It follows from this that $\ga(z) \in C_z\cap C_{z_*},$ and therefore the point $q_*$ belongs to $\cE.$  That $\ga$ is a parametrization of $C_{z_*}$ is clear.

     \bc \label{c1} If  $\Delta(q)\neq 0$ for every $q\in \cD,$ then $\cD= \cE.$ \ec 
   \vs
 \subsection {Conjugate points in Calculus of Variation}  
  
 While the version  $\cE'$ of the envelope is used mainly for the construction of singular solutions to nonlinear differential equations, the version $\cE$ adopted here arises in Calculus of Variations by considering families of extremals through a given point of the plane. 
 Before the introduction of functional spaces in Calculus of Variations, it was customary to identify the functions with the plane curves defined by their graphs. Graphs of function arising as stationary points of the variational functional were termed "extremals". 

Consider  the  variational integral  $ \phi(y) =
\int_0^T  \cL(x, y, y') \, dx, $  with smooth Lagrangian $\cL,$ having   non-vanishing  partial derivative $\tfrac{\partial \cL}{\partial z},$ and assume  that the initial value problem $y(0)=0,\, y'(0) =z$ for the Euler Lagrange equation of the functional $\phi$ has a solution $y(x,z)$  defined for all times. 

   Given a $z_* \in \R$  the graph  $C_{z_*}$ of 
$y_{z_*}(x)=y(x,z_*)$  is naturally embedded  into a one parameter family of extremals  through the point $(0,0)$ defined by  
\[ C_z =\{(x,y) |f(x,y,z):=y-y(x,z)=0\}. \] 
A point  $q_*=(x_*,y_*)\in C_{z_*}$ is said to be \emph{conjugate} to the point $0=(0,0)$ along $C_{z_*}$   if the Jacobi equation (the  linearization of the Euler Lagrange equation of the functional $\phi$ ) has a nontrivial solution vanishing both at $0$ and  $x_*.$ Such a solution  is necessarily of the form $\eta(x)=\frac {\partial y}{\partial z}(x,z_*).$ Hence $ \frac {\partial f}{\partial z} (x_*,y_*,z_*)= -\frac {\partial y}{\partial z}(x_*,z_*)=0,$ and therefore $q_*\in C_{z_*}\cap \cD.$ Conversely, every $q_*\in C_{z_*}\cap \cD$ is a conjugate point.  On the other hand, being $ \frac {\partial y}{\partial z}( - ,z_*)$ a nontrivial solution of a linear differential    equation of second order,   $\frac {\partial ^2 y}{\partial x\partial  z} (x_*,z_*)$ cannot vanish. Since $\Delta( q_*) =\frac {\partial ^2
y}{\partial x\partial  z} (x_*,z_*),$ from Theorem \ref{t2} we get

\bc\label{conextr}
       If  the point $ q_*$ is
conjugate to $0$ along $C_{z_*}$ then for small enough $z$ the extremal $C_{z}$ intersects $C_{z_*}$ at a unique point $ q(z)$ smoothly depending on $z$ and such  $ q(z) \ra  q_*$ as $z \to z_*$. In particular, the conclusion of the Corollary {\rm \ref{c1} } holds in the case of families arising as deformations by extremals. 
\ec

The first conjugate point is of particular importance in Calculus of Variations because at this point the curve loses its minimizing property.  
In  \cite{Ja} Jacobi gave a geometric characterization of the first point conjugate along a given extremal as the point where the extremal touches for the first time the envelope $\mathcal{E}$ of the family. 

 In older texts of calculus of variation the definition of envelope is the one adopted here, see for example the book of Bolza  \cite{Bo}. This also was the original viewpoint of Jacobi \footnote{....Wenn man $l$ an einem Punkt einer Oberfl\"ache nach allen Richtungen k\"urzeste Linien sieht so k\"onnen zwei F\"alle eintreten zwei unendlich nahe k\"urzeste Linien laufen entweder fortw\"ahrend neben einander ohne sich zu schneiden oder sie schneiden sich wiederum und als dann bildet die Continuit\"at aller Durchschnittspunkte ihre einhüllende Curve Im ersten Falle h\"oren die k\"urzesten Linien nie auf k\"urzeste zu sein im zweiten sind sie es nur bis zum Ber\"uhrungspunkte mit der einh\"ullenden Curve...},  but probably because  $\cD$ is considerably easier to deal with, the equality  $\cD= \cE$ in the case of extremals, contributed to blur the distinction between $\cE$ and $\cD$ in general and led many modern authors to refer to $\cD$ as the "envelope". 
 
\br  Families of geodesics through a given point of a Riemannian manifold behave similarly.   Appropriately reformulated,  Corollary \ref{c1} holds in this case as well.   However, on semi-Riemannian manifolds, the envelope of the family of geodesics through a given point can be a proper subset of the discriminant
\cite{MPW, MPP}.
\er
\vs

\subsection{Envelopes, bifurcation, and branching} 
 Our proof of Theorem \ref{t2} can be found in several old analysis books as a typical application of the Implicit Function Theorem.  With the advent of Bourbaki's axiomatization of mathematics, it disappeared along with envelopes and many other shadowy concepts, only to reappear in a slightly different guise in the bifurcation theory of the nineties.
 
As a motivation of our methods in the next section we shall briefly discuss the equivalence between Theorem \ref{t2} and the following  lemma which constitutes the main step in the proof of the well known Crandall-Rabinowitz Theorem in bifurcation theory:
\bl [Bifurcation Lemma]\label{l1}
Let $f\colon I \times \R \ra \R$ be a smooth  function such that $f(\la,0)=0$ for all $\la.$ 

If $\, \frac{\partial f}{\partial x} (\la_*,0)=0$ and
 $ \frac {\partial^2 f}{\partial \la \partial x} (\la_*,0)\neq 0,$ 
 then  there exists  a smooth  function $\la \colon (-\delta, \delta) \ra \R$  such  that  $\la(0) = \la_*$  and $f(\la(x),x) =0.$ \el
 %\vs\vs\vs\vs\vs
% \centerline{\includegraphics[width=0.8\textwidth]{bifurcpic.jpg} }
Here the $\la$-axis  is a "trivial branch" of solutions of the equation $f(\la,x)=0$
 and  $\ga(x) =(\la(x),x)$ is a smooth branch of nontrivial solutions bifurcating the trivial branch at  the point $(\la_*,0).$ 
 
The Envelope Theorem and Lemma \ref{l1} are related by  duality. Namely, if we  exchange the role of $\la$ and $x$ by taking  $x$ as a parameter and $\la,y$ as variables then the Envelope Theorem applied to the family of curves $C_x$ defined by the equation $g(\la,y,x):= y-f(\la,x)=0$ gives the conclusion of Lemma \ref{l1}. Indeed, since $ \frac{\partial g}{\partial x} (\la_*,0)= -\frac{\partial f}{\partial x} (\la_*,0)=0,$ the point $q_*=(\la_*,0)\in C_0$ belongs to the discriminant $\cD$ of the family. Moreover  
  $ \Delta (q_*)=[g_y g_{\la x} - g_\la g_{yx}] ( \la_*,0,0)= \frac {\partial^2 f}{\partial \la \partial x} (\la_*,0)\neq 0,$ and by the Envelope Theorem there exists a curve $\ga (x) = (\la(x), y(x))$  defined in $ (-\delta, \delta)$ with $\ga(0)=q_*$ and such that $\ga(x) \in C_x\cap C_0.$ But $C_0$ is the $\la$ axis and hence 
  $y(x)\equiv 0.$  Thus $f(\la(x),x) =0,$ which proves the lemma.
  
   To reverse the order it is enough to notice that any family of regular curves $f(x,y,z)=0$  by a local change of coordinates can be turned into a family of graphs  $ y- f(z,x)=0.$ To prove  the Envelope Theorem, in this case,  it is enough to apply  the Bifurcation  Lemma 
to $\tilde f (z,x) = f(z,x)- f(z,0),$  with the $z$-axis as the trivial branch. 

\br The bifurcation Lemma can be easily proved along the lines of the proof of the Envelope Theorem.   Simply, write $f$ in the form $f(x,\la)=\la g(x,\la).$ Then $g(x_*,0)=0$ and  $\frac{\partial}{\partial \la} g(x_*,0)  =\frac {\partial^2 f}{\partial x\partial \la} (x_*,0) \neq 0.$ Now apply the Implicit Function Theorem to $g$. \er

The next observation is at the core of our approach to Balasko's equilibrium manifold. 

Let  $C_z = \{(x,y) | f(x,y,z) =0\} $ be an implicitly defined  one-parameter family of plane curves such that  $0$ is a regular value  of the function $f.$  Then $$M= \{(x,y,z) | f(x,y,z) =0\}$$
 is a surface in $\R^3$ and each $C_z$ is the image under projection $\Pi(x,y,z):=(x,y)$ of the intersection of $M$ with the horizontal plane  $H_z$ through $z$.
 
  \bl \label{l2} Let $\pi$ be the restriction of $\Pi$ to $M.$ The discriminant $\cD$ and the envelope $\cE$ of the family $\{ C_z\}, z\in \R$ defined by $f$ is the image by $\pi$ of the set its critical points and branch points respectively. 
\el
  
\proof   Using   $ T_p (M) = \Ker Df (p)$ we have that at  $ p=(x,y,z) \in M$ 

\be \Ker D\pi =\Ker \Pi_{|T_p(M)}= \Ker Df  \cap \Ker \Pi = \Ker DF \ee
where  $F(x,y,z))=(x,y, f(x,y,z)).$ But the Jacobian determinant of $F$ does not vanish at $p$ if and only if $\frac {\partial f}{\partial z}(p)\neq 0.$ Hence  $p$ is a critical point of $\pi$ if and only if  $\pi (p) \in \cD.$

The second assertion follows again by duality:  a point $q_*=(x_*,y_*)$ belongs to $\cE$  if and only  if there exist a sequence of  points $q_n=(x_n,y_n)\in C_{z_n}\cap C_{z'_n}$ with $ z_n\neq z'_n $ but $ z_n- z'_n\ra 0 $ and $q_n, \ra q_*$. Here the parameter is $z$ but if we consider $(x,y)$ as parameter and set $z_*= \lim z_n$ we get  
$p_n=(x_n,y_n,z_n) \neq  (x_n,y_n,z'_n)=p'_n$  with $\pi( p_n) = q_n=\pi(p'_n)$ and  $p_n , p'_n \ra p_*=(x_*,y_*,z_*) \in M$ which means that $p_*$ is a branch point of $\pi$ and $q_*=\pi (p_*).$

%\begin{example}If $E =\{ (x,y,z) | z^3-x-y=0\}$ then $(0,0,0)$ is a critical point of the projection $\pi \dd M\ra \R^2$ but it is not  a branch point being $M$  a graph of a function. Hence  $(0,0)$ belongs to the discriminant $\cD$ of the family but not to the envelope $\cE$. In general, the envelope of a family of curves is only a proper subset of the discriminant.  In the case of families of extremals, both sets coincide due to the particular nature of the Lagrangian.   \end{example}

 Using Lemma \ref{l2}  and the Envelope Theorem we can state a sufficient condition  for a critical point of the projection $\pi$ to be a branch point.  
\bt\label{t11}
Let  $f, M,\pi $ be  as above.  If  $p_*=(x_*,y_*,z_*) $ verifies
   \bit 
 \item [i)] $f_z(p_*)=0.$ 
  \item [ii)] The determinant    
\be \Delta(p_*) = [f_y f_{xz} - f_x f_{yz}](p_*) \neq 0,\ee
\eit 
 then $ p_*\in \cB.$  

 Moreover, there is a  path $\ga_0 (t) = (x(t),y(t),t) \in M$ defined in  $(z_* -\gd, z_*+\gd)$  with  $\ga_0(z_*)= e_*$  and such that also  $\ga_1(t) =(x(t),y(t),z_*)\in M.$ 
 
 \et

\vs

\subsection{The intrinsic derivative}
To formulate our criteria in higher dimensions, we will need the notion of intrinsic derivative of a vector bundle morphism introduced in singularity theory by Porteous \cite{Po}.

Given two smooth vector bundles $E, E'$  over a smooth manifold $M$ and  a vector bundle morphism $\psi \colon E\ra E',$  the \emph{intrinsic derivative} of $\psi$ at a point $x\in M$  in the direction of  a  tangent vector $v\in T_x(M)$ is defined as follows:
 Let $E_x, E'_x$ be the fibers at $x$ of the two vector bundles. On a neighbourhood $U$ of $x$ take two trivializations  $\gt \colon U\times E_x\ra E_{|U}$ and  $\gt' \colon U\times E'_x \ra E'_{|U}$  such that $\gt_x=\Id_{E_x}$ and   $\gt'_x=\Id_{E'_x}.$ 
 Then $\gt'^{-1}\psi\gt$ is a morphism of trivial vector bundles and hence has the form $(x', e) \ra (x', L(x') e )$ where  
$L\colon U \ra \Hom(E_x,E'_x) $ is a smooth map. Let $exp$ be any exponential map defined on a neighbourhood of $0$ in $T_x(M)$ with values in $U,$  Then $l(t)= L(exp(tv)),$  is a smooth path of homomorphisms  defined in a neighbourhood of $0$ in $\R.$ 

\bd The intrinsic derivative  $\partial_v\psi(x) \colon  \Ker \psi_x \ra \Coker\psi_x $  in the direction $v$ is defined by
\be \label{intrinsic} \partial_v\psi(x) = k \frac{dl}{dt}(0)i, \ee
where $i\colon\Ker\psi_x \ra E_x$ 
and $k\colon E'_x \ra \Coker \psi_x$  are the inclusion of the Kernel  and the projection to the Cokernel respectively. 
\ed 

  As its name says, this form of derivation  is independent  of  the choice of  trivializations and the exponential map. Indeed, taking another pair of local trivializations $\tilde \gt $ and  $\tilde \gt' $ of $ E_{|U} \hand E'_{|U}$ the corresponding $\tilde L$ is related to  $L$ by $\tilde L(x') = A(x')L(x')B(x')$ where  $ A(x'), B(x')$  are isomorphisms such that $A(x) =\Id =B(x).$   Taking $ x'=\exp (tv) $ we get $\tilde l(t) = a(t)l(t)b(t)$ with $a(0)=\Id= b(0).$  Hence,  denoting with $ \dot l(t)$ the derivative  $\frac{dl}{dt}(t),$ we have 
\[\dot { \tilde l}(0) = \dot a(0)l(0)b(0) + a(0) \dot l(0)b(0) +  a(0)l(0)\dot b(0)\]
But then  $k \dot { \tilde l}(0) i =  k \dot l (0) i $ since the first term vanish on $\Ker L(x)$ and the  third belongs to $\im L(x).$ The independence from the choice of the exponential map follows easily from the fact that   $D \exp (0_x)  = \Id\dd T_x(M) \ra T_x(M).$ 
  \vs
Given a smooth map $f\colon M\ra N$ between two manifolds, its tangent map or differential  can be regarded as  the vector bundle morphism $Df \colon T(M) \ra f^*T(N)$ between the tangent bundle $T(M)$ of $M$ and the pullback $f^*T(N)$ of the tangent bundle of $N$ by $f.$
\vs
By definition the \emph{intrinsic second derivative} of $f$ at $x\in M$ in the direction of $v\in T_x(M)$ is the intrinsic derivative of the morphism $Df.$ Namely, 

\be \partial^2_v f(x)=  \partial_v Df(x)\ee

\vs

\subsection{The theorem of unavoidable branching} 

 In what follows we  will need few notions related to the local form of  Brouwer's topological degree $\deg(f,U,y)$  defined  on admissible triples (see the Appendix C). If $y$ is a regular value of the map, the Brouwer degree of $f$ on $U$  is given by the formula  \[\deg(f, U, y) = \sum_{i=1}^r \sgn Df(x_i),\]
 where   $ f^{-1} (y) =\{x_1,  \dots , x_r\} $  and $\sgn Df(x_i)=\pm 1,$  depending on whether the isomorphism $Df(x_i) $ preserves or reverses orientation of the tangent spaces.
 
   The \emph{multiplicity} of an isolated  point $x\in f^{-1}(y)$ is defined by  $m(f, x)= \deg(f,V, y),$ where $V$ is any neighbourhood  of $x$  such that $V \cap  f^{-1}(y) =\{x\}.$  In particular, if $x$ is a regular point,  then $m(f, x) = \sgn Df(x)= \pm 1.$ The product formula for the degree implies  that the multiplicity is invariant under composition with a local orientation  preserving diffeomorphism. 
\vs   

Let $M, N$ be smooth  $n$-dimensional orientable manifolds  $f \colon M \ra N$ be a smooth proper map. Let $\cB$ the set of all branch points of $f,$ i.e., points at which f fails to be a local homeomorphism. By the Inverse Function Theorem  $\cB$ is a subset of the set $\cC$ of critical points of the map $f.$

To simplify notations, we will consider the branching of maps from orientable manifolds to Euclidean spaces. This is the only case that we will need here and the general case can be easily recovered from this one. 

The following theorem relating non-degeneracy criteria for critical points with unavoidable branching well illustrates the nature of the topological approach to bifurcation. 
\vs 
\vs

\bt  \label{t4}  Let  $M$ be a connected orientable $m$-manifold and let $f\colon M \ra R^m$ be a smooth proper  map.  If  $\xx\in M$ verifies
   \bit 
 \item [i)] $\Ker Df (x_*)$  is odd-dimensional 
  \item [ii)] There exists a direction  $v\in T_{x_*}(M),$ such that   $\partial^2_v f (\xx)$ is an isomorphism. 
 \eit
  Then $x_*$ is a branch point. Moreover, in every neighbourhood  of $x_*$ there are two regular points  of $x_+\hand x_-\hof f$  belonging to different path components of $M\setminus \cB.$  Equivalently, every path in $M$ from $x_+$ to $x_-$ must intersect $\cB.$ \et

 \br  The regular points $ x_\pm$ have opposite multiplicity.\er
  %\centerline{\includegraphics[width=0.6\textwidth]{unavoid.jpg} }
 
  \proof 
  We first will show that $\xx\in\cB.$
  The idea is to work in coordinates of $M$  at $x_*$ defined by the exponential map.   $\exp \colon U\subset T_\xx(M) \to V\subset M.$  Let us take $I =[-\gd, +\gd]$ with  $\gd>0 $ and $U$  small enough such that,  putting   \be \label{e10} g(t,u)=f(\exp(tv)+u)-f(\exp(tv)),\ee the map  $g\colon I\times U \ra \R^n$  i s well defined. 
  
 Let $L(x)= Df(x) \hand  l(t)=L(\exp(tv)).$ Then $D_u g (0,0)=l(0)=L(\xx).$
  
   Let us choose an oriented basis of  $T_\xx(M)$ of the form 
  $\{ e_1,\dots, e_k,e_{k+1}\dots, e_m\}$  with 
  $  \{e_1,\dots,e_k \}$ a basis of $\Ker L(\xx) $ and compute the determinant  of the matrix of $l(t)$ with respect to this basis and the canonical basis of $\R^m.$
  
    In terms of first order Taylor  expansion of the columns  we have  \cite {Fi} 
 \be \label{e11} \det l(t) = d\cdot t^k + o(t^k),\ee
 where \[d =\det \Big( \dot l(0)e_1, \dots   \dot l(0)e_k,  l(0)e_{k+1},\dots, l(0)e_{m}\Big).\]
 
 But $\dot l(0)$  is injective on $\Ker l(0)$ because $k\dot l(0)$ is injective there. Moreover l(0) is injective on the subspace generated by $ \{e_{k+1},\dots, e_{m}\}.$ By assumption $ii),$ no linear combination of the first $k$-vectors belongs to $\im l(0).$ 
It follows then that the column vectors in \eqref{e11} are linearly independent and hence $d \neq 0.$ Therefore $d(t)$ vanishes only at $0,$  and since $k$ is odd $\det l(\pm \delta )$  have opposite signs. By Krasnoselskij Principle (Appendix C) the map $g$ has a bifurcation point  from the trivial branch $u=0$ in  $I $ which, again by \eqref{e11}, must be $t=0.$  But, if  $(t_n,u_n)\ra (0,0)$ in $I\times U,$   $g(t_n,u_n)=0 $ with $u_n \neq 0,$ then, taking  $x_n = exp(t_nv + u_n) ,  x'_n = t_nv,$ we have  $x_n,  x'_n\ra \xx,$ 
$x_n \neq x'_n$ but $f(x_n) = f(x'_n).$   Therefore $\xx \in \cB.$

In order to prove the second assertion let us consider  $x_\pm= exp(\gd_\pm v).$ By the previous step $x_\pm $ are regular points of $f$ of opposite multiplicity.  They belong to different path components of $M\setminus \cB$  as a consequence of the following  lemma: 

\bl \label{l7} Let   $f\colon M \ra R^m$ be a smooth  proper map from a  connected orientable $m$-manifold $M$ into $R^m.$ If  $x_\pm$  are isolated  in $ f\inv(y_\pm)$  and $m(f, x_+)\neq m(f, x_-),$  then  $x_+ \hand x_-$ belong to different components of $M\setminus \cB.$  \el

\begin{proof} 

Let  $\ga\colon[0,1]\ra M $ be a path with endpoints $x_\pm=\ga(i);\, i=0,1.$ Without loss of generality we can assume that $\ga\dd I \ra M$ is an embedding.  Take an open set $\Omega$ with smooth boundary  such that $\ga(I) \subset \bar \gO,$  $\ga( I) \cap \partial\gO=\{x_\pm\}.$ Then $\ga(I)$ is a neat submanifold of $\bar\gO$ and therefore has a tubular neighbourhood in $\bar\gO$  \cite{H}.  This extends the embedding $\ga$ to an embedding  of the  normal disk bundle \[D =\{ (t,u) \in \ga^* (T(M))| u \perp \dot \ga(t) \hand  \| u\| < \ve \}\] of $\ga(I)$ as its  neighbourhood in $\bar\gO.$   Since $I$ is contractible $D$ is a trivial bundle and taking into account the trivialization we obtain an orientation  preserving diffeomorphism $\psi \colon I  \times D(0,\ve)\subset I\times \R^m \ra M $  such that $\psi(t,0) =\ga(t).$
 Much as before, we define $h\colon I\times D(0,\ve) \ra \R^m$ by   
 \be \label{e120} h(t,u)=f(\psi (t,u))-f(\ga(t)).\ee Clearly $0$ is an isolated zero of both endpoint maps $h_i=h(i,-) , i=0,1 $. Moreover, since the degree is invariant under orientation preserving diffeomorphisms,  $m( h_i,0) = m(f, x_\pm), i=0,1.$ 
By hypothesis we have $m( h_0,0)\neq m( h_1,0)$  and hence  by  Krasnoselskij Principle  $h$ must have a bifurcation point $\ttt$ from the trivial branch. It remains to show  that  $\yy = \ga(\ttt)$ is a branch point of $f$ and hence $x_+\hand x_-$ belong to different path components of $M\setminus \cB.$  Indeed, if  $(t_n,u_n)$  is a sequence of  nontrivial solutions   $h(t_n,u_n)=0, u_n \neq 0 $ converging to $(\ttt, 0),$ then   $  x_n =\psi( t_n, u_n) \neq  x'_n = \ga(t_n),$ but $f(x_n) = f(x'_n) $ and  $x_n, x'_n\ra \xx.$ 
  \end{proof} 
  
\bc \label{c6} With the hypothesis of Theorem {\rm\ref{t4}}, assuming that $y_\pm$ are regular values of $f,$  the following alternative holds: either $ y_\pm= f(x_\pm)$  belong to different path components of $ Y= \R^m\setminus f(\cB)$ or,  for any path $\ga$ in $Y$ from $y_+$ to $y_-$  the  unique lifting $\tilde \ga \colon [0,1] \ra M$ such that $f\circ\tilde \ga = \ga$
and $\tilde \ga(0) = x_+$ must have the second end point  $\tilde \ga(1) \neq x_-.$ 
 \ec
 
\begin{proof}
 In the first case every path between $\pm y$ must intersect $\cE=f(\cB).$ In the second, 
  if $C$ is the component containing both points,  we have that $f\colon f\inv(C) \ra C$ is covering map being a proper surjective local homeomorphism \cite{Chu,Ca}. 
   Therefore $f$ has the unique path lifting property, but by Theorem  \ref{t4} the lifted path cannot have $ x_-$ as the second endpoint. 
\end{proof}

 \vs
\subsection{Manifolds with a nowhere dense set of avoidable branch points}  
  
It is well known that isolated critical points of smooth maps $ f \colon M \ra \R^m$ cannot be branch points if $m\geq 3$ \cite{Pl}.  

The next theorem will explore the kind of restrictions imposed on a map $f$ such that $M\setminus  \cB$ is dense and connected.

The following theorem is a simplified version of Theorem 1.8 in \cite{FP}.

\begin{theorem}\label{t8}Let $M$ be as above and let $f\colon M\ra \R^m$ be a smooth  proper map  such that  $\deg(f)=1.$ Assume that  $M \setminus \cB$ is connected  and dense in $M.$  Then, for every $y\in \R^m,$ $f\inv (y)$ is a \emph{continuum}, i.e., a compact connected set. In particular, if $f\inv (y)$ contains some regular point $x$, then $x$ is the only element of this set. 
 \end{theorem}
 
\begin{proof} By properness $f\inv (y)$ is compact. We will show that it is connected by contradiction. Suppose that $f\inv (y)$ is disconnected. 
Take a cover of $f\inv (y)$ by  two disjoint open sets of $M$, $W_i, i=1,2$ with  $ f\inv (y)  \cap W_i \neq \emptyset , i=1,2.$ 

By the additivity and excision property of the degree
\begin{equation} \label{e110} 1 =\deg(f,W_1,y) +\deg(f,W_2,y). \ee
At  least one of the right hand side members  does not vanishes. Assume  that
\begin{equation} \label{e15} \deg(f,W_1,y) \neq 0.\end{equation}
  By the density hypothesis we can find a point  $\xx\in W_2$ not belonging to $ \cB.$ 

Then, since  $f $  is a local homeomorphism at $\xx,$ we can choose two open neighbourhoods  $U\subseteq W_2$ and $V\subseteq\R^m$ of $\xx$ and $f(\xx)$ respectively such that  $f_{|U}$ is a homeomorphism of $U$ onto $V.$  
 By Sard's Theorem $V$ contains a  regular value of  $y$ of  $f. $ Let $f\inv(y)=\{x_1\dots,x_n\}.$  By the definition of topological degree,
 \begin{equation} \label{e300} 1=\deg_b(f,M,y)=\sum_ {I=1}^n m(f,x_i)\end{equation}
 Since $ M \setminus \mathcal B(f)$ is connected, by Lemma \ref{l7} all the multiplicities $m(f,x_i) $ coincide. Substituting  in \eqref{e300}  we get
$ 1= n\cdot m(f,x_1).$ 
It follows then that $n=1$ and, being  $f_{|U}\dd U\ra V$ a homeomorphism,  $U$ must contain the only solution of  the equation $f(x)=y.$  Therefore  $\deg_b(f,W_1,0) = 0$ in \eqref{e101}, which  contradicts  \eqref{e15}.
  \end{proof}
\vs
\subsection{A reduction property of the intrinsic derivative} 

The reduction property that we  need is as follows:
\bl\label{l8} Let $E, E'$  be two vector bundles over a smooth manifold $M$  and let   $\psi \colon E\ra E'$  be a vector bundle morphism. If  $F'$ is a sub-bundle of $E'$ verifying  \be \label{e12} \im\psi  + F' = E',\ee 
  then $F=\psi\inv (F')$ is a sub-bundle of $E$ and, if  $\psi' \colon F\ra F'$ is the restriction of $\psi $ to $F,$ then for every $v\in T_x(M),$ \be \partial_v\psi(x)= j'\partial_v\psi'(x), \ee where $j'\dd \Coker \psi'_x \ra
  \Coker \psi_x$ is the canonical map $F'_x /\im \psi_x \cap F'_x \cong E'_x/\im\psi_x$ induced by the inclusion.   
\el 
\begin{proof} 
 The fibers of $F=\psi\inv (F')$  are the kernels of the morphism  $\pi \psi \dd E \ra E'/F'$ which is surjective by \eqref{e12}. Since the kernels of a vector bundle epimorphism form a vector bundle,  $F$ is a sub-bundle of $E.$ Since the problem is local and the intrinsic derivative is independent of the trivialization, we will consider  al bundles to be trivialized over a neighbourhood $U$ of $x.$  We put $E= U\times W$ and $F= U\times V$ with $V$ a subspace of  $W,$  and similarly with $E', F'.$ We will denote with $i, i'$ the respective inclusions.  Then $\psi $ on a fiber  at $m \in U$ has the form
 $\psi_m w =(m,L_m w),$ with $L\colon U \ra \Hom(W,W')$ smooth, and we are  denoting  with $L_m$  instead of $L(m)$ the value of $L$ at $m\in U.$ Being  $F=\psi\inv (F')$, it follows that   $L\inv_m (V')=V$  for every $m\in U.$ 
 
   Let  $L'_m\in \Hom(V,V') $ be the restriction of $L_m$ to $V.$ We have  $\Ker L_x = \Ker L'_x,$ and since $\im L'_x= V'\cap \im L_x$ the projections to quotient $ k \dd W' \ra \Coker L_x$ and  $k' \dd V' \ra \Coker L'_x$ are related by   $j'k'= ki'.$ Finally, denoting with $j$ the inclusion of $\Ker L'_x=\Ker L_x $ into $V$  and putting $l(t)=L(exp(tv)),$  
  $l'(t)=L'(exp(tv)),$ we have that  $\dot l(0)i=i' \dot l'(0)$ and therefore, 
 
\be \partial_v \psi (x) = k\dot l(0)i j =ki' \dot l'(0) j=j'k'\dot l'(0)j=j'\partial_v \psi' (x).\ee  
\end{proof}

We will use the above reduction property in the next subsection  together  with the following invariance  property  of the intrinsic derivative which is a consequence of the chain rule. \bl\label{l9} Let $E, E'$  be two vector bundles over a smooth manifold $M$  and let   $\psi \colon E\ra E'$  be a vector bundle morphism. Let $N$ be  a submanifold of $M,$ and $\psi_{|N}\dd E_{|N} \ra E'_{|N}$ be the restriction of $\psi$ to $N,$ then, for every $n\in N$ and $v\in T_n(N),$
\be \partial_v \psi =  \partial_v (\psi_{|N}).\ee
    \el
We leave the proof  to the reader.

\vs

\subsection{Intrinsic derivative on implicitly defined manifolds}
\bd An \emph{implicitly defined family of varieties} in $\Rm$  parametrized by $\Rp$ is a family  of the form $C_p=\{ x\in \Rm | f(x,p)=0\} , p \in \Rp$  where  $f\dd \Rmp\ra\Rp$ a smooth map.
\ed 
 In what follows we will always assume  that $0$ is a regular value of $f,$ i.e., the Frechet differential $Df(x,p)$ is surjective for every $(x,p)\in f\inv (0).$  By the Implicit Function Theorem, $M=  f\inv (0)$ is an $m$-dimensional submanifold of $\Rmp$  which is orientable because  $Df$ induces a trivialization of its  normal bundle. Moreover $C_p$ is the image under the projection map $\Pi \dd\Rmp\ra \Rm$ of the intersection  $ (\Rm\times\{p\})\cap M.$ 

Let us denote by $\pi $ the restriction of $\Pi$ to $M.$ Guided by the case of families of plane curves (m=2, p=1) we define the \emph{discriminant} $\cD$ of the family as the set of critical values of $\pi$ and its\emph{ envelope} $\cE$ as the image by $\pi$ of the set $\cB$ of its branch points. We leave to the reader the reinterpretation of the envelope in terms of the family $C_p, p\in \Rp$ analogous to the case of plane curves.

  A comprehensive study of envelopes along the lines of Thom's paper \cite{Th} would require a stronger regularity condition, namely that   $D_x f(x,p)$  is surjective at points of  $M,$ which, together with properness, would then imply that each $C_p$  is a submanifold of $\Rm$ and $\cup_p C_p$ is a fiber bundle.  But here this is not needed since the only goal  of this subsection is to compute the intrinsic second derivative $\partial_v^2\pi(\mst)$  at  a critical point $m_* =(x_*,p_*)$ of $\pi,$ in terms of the map $f.$ 

\vs

 For this,  let  us consider  the map 
$ g=(\Pi,f) \dd\Rmp \ra \Rm\times\Rp,$  
and its differential  $Dg \dd T(\Rmp) \ra g^*[T(\Rm\times\Rp)].$ 
Let $E,E'$ be the restrictions to $M$ of the trivial bundles $T(\Rmp)$  and $g^*[T(\Rm\times\Rp)]$ respectively. 
 We claim that  $Dg_{|M}$ verifies the condition \ref{e12} with respect to the sub bundle  $F'= M\times \Rm\times \{0\}$ of $E'.$  Indeed,  given a point $m\in M$  and an element $ e'=(f',w')\in E'_m $ by surjectivity of $Df(m),$ we can find an $e\in E$ such that  $Df(m)e= w'. $ But then $e'$ can  be written as  $e' = (\Pi e, Df(m)e) + ( f'-\Pi e,0),$ which proves \ref{e12}. 
 \vs
 
  Denoting with  $\pi$ the restriction of the projection $\Pi $ to $M$ we have  $D\pi= \Pi_{|TM}.$  Therefore $Dg_{|M} $ coincides with $D\pi$ as a morphism from  $F=Dg_{|M}\inv(F')= \Ker Df = T(M)$ to $F'.$  
By Lemma \ref{l9} and \ref{l8} for every $v\neq 0 $ belonging to $T_{\mst}(M)$  we have 

\be \label{e100}   \partial_v Dg (\mst)= \partial_v Dg_{|M}(\mst)= j_1\partial_vD\pi(\mst),\ee

where $j_1$ is the canonical isomorphism of $\Coker D\pi(\mst)$ with $\Coker Dg(\mst).$
\vs
On the other hand  $\xst$ is  a regular value of the projection $\Pi\dd \Rmp\ra\Rm,$  with  $\Pi\inv (\xst) = \{\xst\}\times \Rp,$ and the restriction of $Dg$ to $\{\xst\}\times \Rp:= N$ tautologically verifies  the condition \ref{e12} with respect to the sub bundle $F'= N\times \{0\}\times\R^p.$ Moreover   $  F= Dg_{|N}\inv(F') = N\times\{0\}\times\R^p, $ and  hence   $Dg_{|N} \dd F\ra F'$ coincides with $Df_{|N}.$  
Applying  Lemmas  \ref{l9} and \ref{l8} again,  for every $v\neq 0 $ belonging to $T_{\mst}(N)$  we have 

 \be \label{e101}   \partial_v Dg (\mst)= \partial_v Dg_{|N}(\mst)= j_2 \partial_vDf_{|N}(\mst).\ee

where $j_2$ is the canonical isomorphism of $\Coker Df_{\xst}$ with $\Coker Dg(\mst).$
\vs
Since  $\Ker D\pi(\mst)= \Ker \Pi_{|T_{\mst}(M)}=\Ker \Pi \cap T_{\mst}(M),$ every vector $v \in\Ker D\pi(\mst)$ is tangent to both $M$ and $N.$  Therefore,  putting  $ j' =j_1\inv j_2,$ from  identities  \eqref{e100} and \eqref{e101} we get
\be \label{e102}    \partial^2_v\pi(\mst) =\partial_vD\pi(\mst) = j'\partial_vDf_{|N}(\mst)= j'\partial^2_vf_{|N},\ee for every $ v\in \Ker D\pi.$

Let  $f_{\xst}\dd \Rp\ra \Rp$ be the map $f_{\xst}(p):=f(\xst,p).$ Identifying $N$ with $\Rp$ via the map $j(p) =(\xst,p),$ the morphism   $Df_{|N}$ coincides with the morphism  $Df_{\xst}.$
With the same identification 
 \be \label{e103}  \Ker D\pi(\mst) =\Ker \Pi \cap T_{\mst}(M) =\{0\}\times\Rp \cap \Ker Df(\mst)= \Ker Df_{\xst}(\pst).\ee 
Clearly $\Coker Df_{|N}= \Coker Df_{\xst},$ and hence we have
 
\bt \label {t20} A point $ \mst=(\xst,\pst)\in M$ is a critical point of $\pi$ if and only if $\pst$ is a critical point of $f_{\xst}.$ Morever, $\Ker D\pi(\mst)= \Ker Df_{\xst}(\pst)$ and  for every $ v\in \Ker Df_{\xst}(\pst)$
\be \partial ^2_v\pi(\mst) =j'\partial^2_vf_{\xst}(\pst) \ee \et

\vs
\subsection{A criterion for unavoidable branching in terms of the defining map. }

\bt\label{t21}
 Let  is $f\colon \Rm\times\Rp \ra R^p$ be a smooth map  such that $0$ is a regular value of $f.$  Assume that  $M=f\inv(0) $ is  connected and locally bounded over $\Rm.$  
 
  If  $\mst = (\xst, \pst) \in M$ verifies
   \bit 
 \item [i)] $\Ker Df_{\xst} (\pst)$  is odd-dimensional. 
  \item [ii)]  There exists a vector  $v\in \Ker Df_{\xst} (\pst),$ such that 
 \[D^2f_{\xst}(\pst)[v,u] \notin \im  Df_{\xst}(\pst) \hforall u \in \Ker Df_{\xst}(\pst), u\neq 0.\] 
 \eit
  Then $m_*\in \cB.$ Moreover, in every neighbourhood  of $m_*$ there are two regular points  of $m_+\hand m_-\hof f$  belonging to different path components of $M\setminus \cB.$  In particular, $\cB$ separates $M.$ 
  \et
  \begin{proof} That $M$ is locally bounded over $\R^m$ means that every point $y\in \Rm$  has a neighbourhood $U$ such that $M\cap U\times \Rp$ is bounded in $\Rmp.$ It follows then that the projection $\pi\dd M\ra \Rm$ is a proper map.
Checking the  definition of the intrinsic derivative  we find that  for every $u,v\in\Ker Df_{\xst}(\pst)$ 
\be  \partial^2_vf_{\xst}(\pst)u = q D^2f_{\xst}(\pst)[v,u],\ee  
where $D^2f_{\xst}(\pst)$ is the second differential  of $f_{\xst}$ at $\pst$ and $ q\dd\Rp\ra \Coker Df_{\xst}(\pst)$ is the projection to the quotient. By dimensional reasons 
$q D^2f_{\xst}(\pst)[v,-]$ is an isomorphism if and only if it is a monomorphism, or equivalently if  
\be D^2f_{\xst}(\pst)[v,u] \notin \im  Df_{\xst}(\pst) \hforall u\in \Ker Df_{\xst}(\pst) \ee 
Now the  claimed result follows from Theorems \ref{t20} and \ref{t4}
\end{proof}
\vs
 \section{APPENDIX}
 \vs
 \subsection{A: Equilibrium Prices and Crisis}
 We cannot resist but to quote word by word the vivid description of the relation be between the discontinuity of equilibrium prices and an imaginary crisis of the market of artichokes in Bretagne by Balasko \cite{Ba88}.
   
  \includegraphics[width=0.8\textwidth]{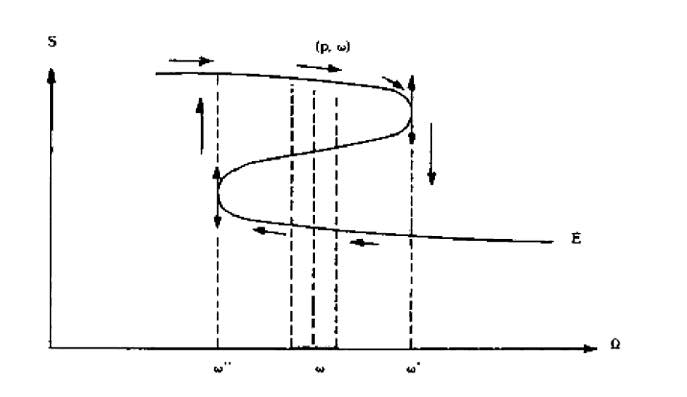}
  
\emph{
"Illustrons cette importance par un exemple quelque peu imaginaire mais qu'il n'est pas difficile de relier \`a une r\'ealit\'e plus concr\`ete. Interpr\'etons le march\'e consid\'er\'e comme \'etant celui des artichauts sur un lieu de production comme la Bretagne et examinons la s\'equence des prix d'\'equilibre observ\'es jour apr\`es jour. Au d\'ebut de lasaison, les artichauts sont disponibles en petites quantit\'es repr\'esent\'ees par une valeur petite du parametre $w$ ce qui dans le cas de figure donne un prix d'\'equilibre unique et relativement \'elev\'e. Chaque jour nouveau voit les quantit\'es d'artichauts mises sur le march\'e augmenter, $w$ se d\'eplace done vers la droite. Tant que $w$ n'a pas franchi $w"$ il y a unicit\'e du prix d'\'equilibre. En outre, les variations de prix que l'on constate sont faibles : l'accroissement des quantit\'es mises sur le march\'e ne se traduit pas par une baisse sensible des cours. Le franchissement de correspond \`a l'entr\'ee dans une zone d'\'equilibres multiples : trois \'equilibres sont alors possibles, mais le principe de continuit\'e d\'ej\`a vu permet de comprendre que le prix qui va   \^etre effectivement observ\'e sur le march\'e est celui qui correspond \`a la branche sup\'erieure de la courbe $E.$  Autrement dit. l'observation du syst\`eme de prix ne permet pas de mettre en \'evidence, done d'identifier, une modification quelconque des caract\'eristiques du march\'e. Laissons les jours passer. Voici que maintenant $w$ s'approche de $w'.$ Pour le moment rien de significatif n'est observ\'e si ce n'est des variations un peu plus fortes du syst\`eme de prix qui, bien que restant d'un ordre de grandeur relativement petit, n'en sont pas moins nouvelles. Ceci est du au fait que la pente de la tangente \`a $E$ augmente r\'eguli\`erement en valeur absolue. Des observateurs d\'ecrivant cette phase feraient \'etat d'une certaine nervosit\'e du march\'e. Puis arrive le jour où $w$ d'inf\'erieur \`a $w'$ - lui devient sup\'erieur. Le prix d'\'equilibre ne peut plus \^etre prolong\'e infinit\'esimalement. En effet, pour la nouvelle valeur du param\`etre $w.$ Il n'y a plus qu'un seul prix d'\'equilibre possible et sa valeur est tr\`es diff\'erente de l'ancienne. Quelle situation observe-t-on alors sur le march\'e ?}
\emph{On pourrait imaginer qu'il a suffit d'un artichaut pour franchir $w'$ de sorte qu'il est pratiquement impossible pour un observateur \'economique de diff\'erencier les deux jours en question. Mais le march\'e fait cette diff\'erence quand il constate que les prix de la veille n'\'equilibrent plus l'offre et la demande. Le march\'e commence par chercher \`a r\'etablir l'\'equilibre en proc\'edant \`a de petites variations des cours. Mais contrairement aux jours pr\'ec\'edents, aucune des petites variations possibles ne permet d'aboutir \`a un \'equilibre de l'offre et de la demande. Par cons\'equent, apr\`es un certain temps perdu dans ces efforts infructueux (il arrive souvent que le fonctionnement du march\'e soit alors suspendu afin de r\'etablir une certaine s\'er\'enit\'e dans les esprits face \`a cette incapacit\'e du march\'e \`a s'\'equilibrer), la constatation s'impose que si \'equilibre il y a, il va être tr\`es diff\'erent de celui de la veille. Pour les observateurs qui appliquent le principe que les m\^emes causes ont les m\^emes effets, cette discontinuit\'e est inexplicable. C'est donc dans une atmosph\`ere faite d'incompr\'ehension voire d'incr\'edulit\'e face aux \'ev\'enements que l'on voit le march\'e s'ajuster sur un nouvel \'equilibre. Cette phase est souvent accompagn\'ee par une remise en cause passionnelle du m\'ecanisme de marche, remise en cause justifi\'ee par cette discontinuit\'e en apparence impr\'evisible et inexplicable. L'aspect passionnel peut aller jusqu'\`a une destruction de ressources dans le but de faire remonter les cours : qui n'a entendu parler d'artichauts mis \`a la d\'echarge publique ?"}

\vs

 \subsection{B: Examples of Critical Economies} 
  
By Sonnenschein-Mantel-Debreu Theorem every continuous map verifying on $S$ the Walras law coincides on any compact subset of $S$  with the aggregate excess demand of a market with $m\geq l,$ and fixed endowments. 

In order to find a critical economy that admits a continuous selection of prices, it suffices to consider a  market with two participants and two commodities.  Normalize the price of the second good to $p_2=1$ and take a reduced excess of demand, i.e., the excess of demand for the first commodity,  of the form $\bar z(p,\omega ) =(p-1)^3 -(\go_2 - 1 ),$  where $\go= (\go_1, \go_2)$ is the aggregate endowment. Its equilibrium manifold $E$  is a graph. Each  point   $e=(p,\go_1,1,)$ with $p=1$ is a critical equilibrium which not a crisis, since  $E$ posses a global continuous selection of prices $p(\go) =  1+ \sqrt[3]{\go_2 -1}.$  There is a delicate point in this reasoning because we are looking at  $z$ as a function of both prices and endowment. However, since the problem is local, the version of the Sonnenschein-Mantel-Debreu Theorem proved in \cite[Theorem 4.1]{CE} suffices to guarantee the existence of utility functions whose reduced excess of demand is the above map. 

 Several examples of critical economies with continuous price selections can be constructed on markets with two goods using quasilinear utilities of a special form considered in  \cite{ MGW, APP}. 
For example in Figure 4 of  \cite{APP} is represented such an equilibrium manifold $E$ separating regions of the uniqueness of equilibria from the region filled with manifolds with multiple equilibria.  Much like the one described above this manifold is a graph of a continuous function having a vertical tangent at a critical economy. The authors           
 of \cite{APP} used the Negishi approach which substitutes prices with social welfare weights but the results hold for prices as well due to their strict relation with the social welfare weights. 
 
The next example illustrates the criterion of unavoidable branching stated in Theorem \ref{t21}. The example comes from \cite[page 521]{MGW} taken as an illustration of the existence of multiple equilibria.
The two  consumers have their preferences represented  by a quasilinear  utility  functions of  type
 \be \begin{cases} u_1(x,y) = x-  \frac{1}{\ag}y^{-\ag}  \\ u_2 (x,y) = y- \frac{1}{\ag}x^{-\ag} \end{cases} \ee
 Let   $ (\go_i,\go'_{i}),  i=1,2 $  be the  initial endowment of the consumer $i.$
 Choosing the second good as the  numeraire the wealth of the consumer $i$ is  $w_i=p\cdot \go_i+1\cdot \go'_i.$ 

The demand of each consumer is easily obtained using Lagrange multipliers for the constrained maximization of  $u_i$ on the budget set  $px+ 1y =w_i.$  For the reduced excess demand we only need to consider the first good. 

 We get 
\be \bar z(p,\go)= p\inv \go'_1  - p^{-\frac{\ag}{\ag+1}}+ p^{-\frac{1}{\ag+1}} - \go_2 \ee 
With $\ag =8$ and taking   $\go'_1 = \go_2= 0.766$  it turn out that  $\bar z(p)$ vanishes together with its derivative  at approximately  $0.54$ but  $z''(0.54)\neq 0.$ therefore any endowment verifying  $\go'_1 = \go_2= 0.766$ is an unavoidable  crisis. 

 \includegraphics[width=0.6\textwidth]{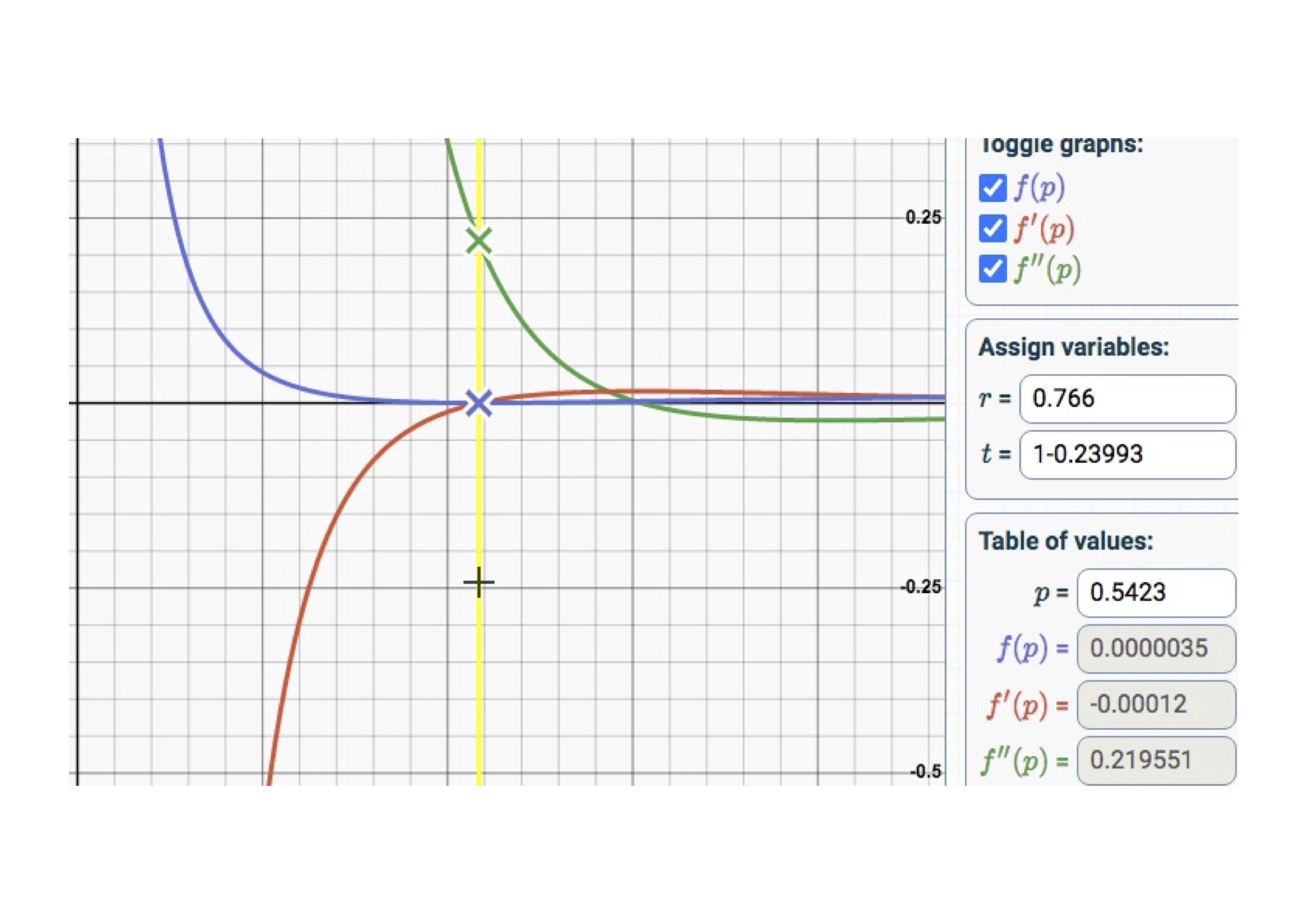}

\vs 
\subsection{C: Brouwer's Degree and Krasnoselskij Principle}  We will use a local version of topological degree on oriented smooth manifolds.  
Given two oriented $n$-manifolds $M \hand N,$   an admissible triple  $(f,U,y)$ is a triple defined  by three admissible elements. Namely, an open subset $U$ of $M,$ a continuous map $f \colon U \ra N$  and $y\in N$ such that $f\inv(y)$ is a compact subset of $U.$  

\begin{theorem} There is one and only one function $\deg$ which assign to every admissible triple an integer $deg (f,U,y)$ verifying: 
\begin{itemize} 	 
 \item[\rm Existence:]\   If   $\deg(f,U,y) \neq 0,$ the equation  $f(x)=y $ has a solution in  $U$. \\
 \item[\rm Homotopy:]\  If $h \colon [0,1] \times U \ra \R^n $ is a homotopy such that   $h\inv(y)$ is compact, then  $$\deg (h_0,U,y)=\deg (h_1,U,y).$$
 \item[\rm Additivity:]  If  $U_i, i=1,2$  are  disjoint  admissible sets, then 
 $$\deg (f,U_1 \cup U_2,y)=\deg (f,U_1,y)+\deg (f,U_2,y),$$ 
  \item[\rm Excision:]\ If  $U\subset V,$ are open sets and   $f\inv(y)$ is a compact subset of $U,$ then $$\deg (f,V,y)=\deg (f,U,y).$$
  \end{itemize}
  \end{theorem}
  The usual  construction is made  by assuming first  that  $f $ is smooth and $y$ is a regular value of $f$ in $U.$ In this case  the degree is explicitly defined  by the formula 
 \be \label{degree} \deg(f,U,y) = \sum_{i=1}^r \sgn Df(x_i),\ee
          where   $ f^{-1} (y) =\{x_1,  \dots , x_r\} $  and $\sgn Df(x_i)=\pm 1$ depending on whether the isomorphism $Df(x_i)$  preserves or reverses orientation of the tangent spaces.  
          The extension  to  the  general case is obtained through approximation of continuous maps by smooth maps and the use of  Sard's Theorem to approximate $y$ by regular values of the map $f.$ 
          
    If  $f\colon M\ra N$  proper and N is connected,  then, by homotopy invariance, the local degree $\deg (f, M, y)$ is independent of the choice of $y\in N$  and will be denoted with $\deg (f).$      

\sk

  Let $I \subset \R $ be an interval and let  $f\colon I\times \R^n\rightarrow \R^n$  be a continuous family of $C^1$ maps  such that $f(\la,0)=0 \hforall \la\in  [a,b].$ the subset $ \cT = I\times\{0\}$  of the set of all solutions of the equation $f(\la,x)= $ is  called the  {\it trivial branch}. 
We will denote with $f_\la \colon \R^n\ra \R^n$ the map defined by $f_\la(x)=f(\la,x).$
 In bifurcation theory we look for solutions of the equation $f(\la,x)=0$ that are arbitrarily close, but do not belong to the trivial branch.
 \bd A {\it bifurcation point  from the trivial branch} of  solutions of  the equation $f(\la,x)=0$ 
is a point $\la_*$ in $I$  such that every neighbourhood of $(\la_*,0)\in \LL\times \R^n$ contains a nontrivial solution $(\la,x) , x\neq 0,$ of  this equation. \ed

The following lemma is an easy consequence of the homotopy invariance of topological degree 
 \bl [Krasnoselkij Principle]\label{l4}
 Let $ f\colon [0,1] \times \R^n \ra \R^n$  be a family of  $C^1$ compact maps such that  $f(\la,0) =0.$
If  both $Df_0(0),\, Df_1(0) $ are invertible  and  \[\sgn\det Df_0(0) \neq \sgn\det Df_1(0),\] then  the interval $(0,1)$  contains at least one  bifurcation  from the trivial branch.  \el

 {\bf Proof } If there are no bifurcation points in $(0,1),$ for small enough $r,$ the only solutions of $f(\la,x)=0$   in $[0,1]\times \bar B(0,r)$ are the trivial ones. But then, by homotopy invariance   $\deg (f_{i},B(0,r),0) =  \sgn\det Df_i(0) ,\, i=0,1,$ b contradicting the hypothesis.

  The extension of the above lemma to the maps between oriented vector spaces needed in the proof of Theorem \ref{t4}  is obvious. 
  
  \vs

\vs

 \end{document}